\documentclass[12pt,amssymb]{article}

\usepackage{latexsym}
\usepackage{amsmath}
\usepackage{amsfonts}
\usepackage{amssymb}
\usepackage{amsthm}
\usepackage{amscd}

\input xy
\xyoption{all} \hfuzz=1.5pt \tolerance=500

\newcommand{\C}{\mathbb{C}}
\newcommand{\N}{\mathbb{N}}
\newcommand{\R}{\mathbb{R}}

\newcommand{\va}{\varphi}

\newcommand{\id}{{\bf 1}}
\newcommand{\co}{{\B C}}
\newcommand{\s}{\sigma}
\newcommand{\lm}{\lambda}

\newcommand{\bd}{\begin{document}}
\newcommand{\ed}{\end{document}}
\newcommand{\Htp}{\mathbin{\stackrel{\cdot}{\bigotimes}}}
\newcommand{\Hts}{\mathbin{\stackrel{\cdot}{\bigoplus}}}
\newcommand{\nl}{\nu\in\Lambda}
\newcommand{\xl}{X_\nu;\nl}
\newcommand{\x}{X_\nu}
\newcommand{\la}{\langle}
\newcommand{\ra}{\rangle}
\newcommand{\ch}{\cal H}
\newcommand{\su}{\subseteq}
\newcommand{\B}{\Bbb}
\newcommand{\e}{\varepsilon}
\newcommand{\ov}{\overline}
\newcommand{\vk}{\varkappa}
\newcommand{\vt}{\vartheta}
\newcommand{\al}{\alpha}
\newcommand{\de}{\delta}
\newcommand{\lo}{\Longleftrightarrow}
\newcommand{\cK}{{\cal F}}
\newcommand{\bb}{{\cal B}}
\newcommand{\dd}{{\cal D}}
\newcommand{\cR}{{\cal R}}
\newcommand{\f}{{\cal K}}
\newcommand{\cb}{{\cal CB}}
\newcommand{\wrr}{\widetilde{\cal R}}
\newcommand{\br}{^\bullet{\cal R}}
\newcommand{\en}{E_\nu}
\newcommand{\el}{{\cal L}}
\newcommand{\lr}{\Longrightarrow}
\newcommand{\Long}{\Longleftarrow}
\newcommand{\q}{\quad}
\newcommand{\qq}{\qquad}
\newcommand{\cd}{\cdot}
\newcommand{\fur}{f: E\to{\B R}}
\newcommand{\furo}{f_0: E_0\to{\B R}}
\newcommand{\fuc}{f: E\to{\B C}}
\newcommand{\ii}{\infty}
\newcommand{\di}{\diamondsuit}
\newcommand{\bgd}{{\bigtriangledown}}
\newcommand{\bu}{{\bigtriangleup}}
\newcommand{\bc}{{completely bounded}}
\newcommand{\cc}{{completely contractive}}
\newcommand{\qs}{{quantum space}}
\newcommand{\isc}{{isometric}}
\newcommand{\ism}{{isomorphism}}
\newcommand{\qss}{{quantum spaces}}
\newcommand{\bco}{{completely bounded operator}}
\newcommand{\bcos}{{completely bounded operators}}
\newcommand{\res}{{respectively}}
\newcommand{\tp}{{tensor product}}
\newcommand{\eq}{{equivalent}}
\newcommand{\qtp}{{quantum tensor product}}
\newcommand{\mma}{\mathrel{\mathop{\otimes}\limits_{A}}}
\newcommand{\mmb}{\mathrel{\mathop{\otimes}\limits_{\bb}}}
\newcommand{\mmp}{\mathrel{\mathop{\otimes}\limits_{p}}}
\newcommand{\mmh}{\mathrel{\mathop{\otimes}\limits_{h}}}
\newcommand{\mmf}{\mathrel{\mathop{\otimes}\limits_{4}}}
\newcommand{\mmi}{\mathrel{\mathop{\otimes}\limits_{i}}}
\newcommand{\mmm}{\mathrel{\mathop{\otimes}\limits_{\bb-\bb}}}
\newcommand{\mmo}{\mathrel{\mathop{\cdot}\limits_{1}}}
\newcommand{\mmA}{\mathrel{\mathop{\otimes}\limits_{A-A}}}
\newcommand{\mmd}{\mathrel{\mathop{\cdot}\limits_{2}}}
\newcommand{\msp}{\mathrel{\mathop{\otimes}\limits_{sp}}}
\newcommand{\mms}{\stackrel{h}{\otimes}}
\newcommand{\mmt}{\stackrel{4}{\otimes}}
\newcommand{\mmx}{\stackrel{i}{\otimes}}
\newcommand{\mmy}{\stackrel{p}{\otimes}}
\newcommand{\mmz}{\stackrel{sp}{\otimes}}
\newcommand{\gd}{\ddagger}
\newcommand{\od}{\odot}
\newcommand{\mt}{\mapsto}
\newcommand{\mmc}{\mathrel{\mathop{\otimes}\limits_{\sim}}}
\newcommand{\mme}{\stackrel{\sim}{\otimes}}
\newtheorem{thm}{Theorem}[section]

\newtheorem{lem}[thm]{Lemma}
\newtheorem{prop}[thm]{Proposition}
\newtheorem{cor}[thm]{Corollary}

\begin{document}

\centerline{{\bf Metric version of projectivity for normed modules}}

\centerline{{\bf over sequence algebras}\footnote{This research was
supported by the Russian Foundation for Basic Research (grant No. 08-01-00867).}}

   \vspace{1cm}

\centerline{A.~Ya.~Helemskii}

\centerline{Faculty of Mechanics and Mathematics}

\centerline{Moscow State University}

\centerline{Moscow 119992 Russia}

\vspace{1cm}

{\bf Abstract.} We define and study the metric, or extreme version of the notion of a projective normed module.
The relevant definition takes into account the exact value of the norm of the module in question, in contrast with
the standard known definition that is formulated in terms of norm topology.

After some preliminary observations and discussion of the case where the base normed algebra $A$ is just $\C$, we concentrate on the case of the next degree of complication, where $A$ is a sequence algebra, satisfying some natural
conditions. The main results give a full characterization of extremely projective objects within the
subcategory of the category of non-degenerate normed $A$--modules, consisting of the so-called homogeneous modules.
We consider two cases, `non-complete' and `complete', and the respective answers turn out to be essentially different.

In particular, all Banach non-degenerate homogeneous modules, consisting of sequences,  are extremely projective within the category of Banach non-degenerate homogeneous modules. However, neither of themis extremely projective within the category of all normed non-degenerate homogeneous modules, when it is infinite-dimensional. On the other hand, submodules of these modules, consisting of finite sequences, are extremely projective within the latter category.

\bigskip
\centerline{{\bf Introduction: Formulation of the main results; comments}}

\bigskip
The concept of a projective module is one of the most important in algebra. In particular, it plays the role of one of
the three pillars of the whole building of homological algebra; the other two are the notions of injective module and a flat module. The first functional-analytic versions of the three mentioned notions appeared about 40 years ago~\cite{matsb,fl,1book}. They were introduced in connection with the rise of interest, in functional analysis, to such topics as derivations of Banach algebras, their extensions and amenability.  The relevant definitions were given in the framework of a certain kind of relative homology, adapted to the  context of functional analysis.
They were formulated in terms of the norm topology of modules in question rather than the norm itself.

Quite recently, the birth of new areas of analysis, notably of the so-called quantum functional analysis (= operator space theory), has caused the introduction and study of metric, or extreme versions of these notions (cf.~\cite{he4,wit}). The specific feature of the new versions is that they take into account the exact value of the norm.

Before the proceeding to our main definition, we recall and fix some terminology and notation.

Let $A$ be a complex algebra, and $X,Y,Z$ left $A$--modules. When $A$ is fixed,
the term `morphism', will always mean a morphism of $A$--modules. Let $\tau:Y\to X$ be a surjective morphism, and $\va:Z\to X$ an arbitrary morphism. A morphism $\psi:Z\to Y$ is called a {\it lifting of $\va$ with respect to $\tau$}, if the diagram
$$
\xymatrix@R-10pt@C+15pt{
& Y \ar[d]^{\tau}\\
P \ar[ur]^{\psi} \ar[r]^{\va} & X }
$$
is commutative.

From now on we suppose that $A$ is a normed algebra, not necessary unital. We denote by $A-mod$ the category of all left normed $A$--modules and their bounded morphisms. Throughout this paper, all normed algebras and  modules are always supposed to be contractive; this means
that for $A$ and $X\in A-mod$ we always have the multiplicative inequalities $\|ab\|\le \|a\|\|b\|$ and $\|a\cd x\|\le \|a\|\|x\|$ for all $a,b\in A, x\in X$.

In what follows, a left $A$-module $X$ is called {\it essential} (they often say also `non--degenerate'), if the closure of the linear span of the set $\{a\cd x: a\in A, x\in X\}$ is the whole $X$.
A bounded morphism $\s$ of $A$-modules is called {\it near--retraction}, if it is contractive (that is $\|\s\|\le1$) and, for every $\e>0$, it has a right inverse morphism with norm $<1+\e$. A normed $A$-module $X$ is called {\it near--retract} of an $A$-module $Y$, if there exists a near--retraction from $Y$ onto $X$.

As usual, the category of all normed spaces and bounded operators is denoted by
{\bf Nor}, and its full subcategory, consisting of Banach spaces, by {\bf Ban}.

Recall that an operator, in particular, a module morphism, is called {\it coisometric} (they say also `quotient map'), if it maps the open unit ball of the domain space onto the open unit ball of the range space. Equivalently,
$\tau:Y\to X$ is coisometric, if it is contractive, and, for every $x\in X$ and $\e>0$, there exists $y\in Y$ such that $\tau(y)=x$ and $\|y\|<\|x\|+\e$.

Finally, our definition takes into account a certain full subcategory in $A-mod$, so far arbitrary chosen. We denote it by ${\f}$.

\medskip
{\bf Definition I.} A module $P\in A-mod$ is called {\it extremely (or metrically) projective with respect to $\f$},
if, for every coisometric morphism $\tau:Y\to X$ of modules in $\f$, an arbitrary bounded $A$-module morphism $\va:P\to X$, and $\e>0$, there exists a lifting $\psi:P\to Y$ of $\va$ with respect to $\tau$ such that $\|\psi\|<\|\va\|+\e$.

\medskip
{\bf Remark.} It is easy to see that the given definition is equivalent to the following one: $P$ is
extremely projective with respect to $\f$ if for every coisometric morphism $\tau:Y\to X$ in $\f$,  
the operator ${\bf h}_A(P,\tau):{\bf h}_A(P,Y)\to{\bf h}_A(P,X)$, acting by the rule $\va\mt\tau\va$, is 
also coisometric. (Here ${\bf h}_A(\cd,\cd)$ is the symbol of the space of bounded morphisms, equipped 
with the operator norm). Thus, in the concise categorical language, $P$ is  extremely projective with 
respect to $\f$, 
if the covariant morphism functor  ${\bf h}_A(P,?):\f\to{\bf Nor}$ preserves coisometries.

\medskip
Extremely projective modules are the principal objects of this paper. Nevertheless, sometimes we shall come across another, somewhat weaker property than that of extreme projectivity. Namely, a module $P\in A-mod$ is called {\it topologically projective with respect to $\f$} if for every $\tau:Y\to X$ and $\va:P\to X$ as before, there exists a (just) bounded lifting of $\va$ with respect to $\tau$.

We would like to emphasize that we assume our modules to be Banach (= complete) only if it is explicitly stated.
(Indeed, we shall see that some of the results and constructions in this paper deal with essentially non-complete modules).

\medskip
First, let us look at the simplest case $A={\C}$, where our modules are just normed spaces.
In this situation, if a normed space $P$ is extremely projective with respect to {\bf Nor} (that is, to the whole ${\C}-mod$), we just say that $P$ is {\it extremely projective in {\bf Nor}}.

Extremely projective normed spaces can be characterized by the following
statement. For an arbitrary non--empty set $\Lambda$, we denote by $l_1^0(\Lambda)$ the
space of all finitely supported functions on $\Lambda$, equipped with the $l_1$--norm. (Thus $l_1^0(\Lambda)$ is a dense subspace in the standard Banach space $l_1(\Lambda)$.) Besides, we set $l_1^0(\emptyset):=l_1(\emptyset):=0$.

\medskip
{\bf Proposition}. {\it A normed space is extremely projective in {\bf Nor} if and only if it is a 
near-retract in {\bf Nor} (= $\co$-mod) of $l_1^0(\Lambda)$ for some index set $\Lambda$.}

\medskip
If we are given a Banach space, we say that it is {\it extremely projective in {\bf Ban}}, if it is extremely projective with respect to {\bf Ban}. These spaces have much more transparent description (cf.~\cite{gro}):

\medskip
{\bf The Grothendieck Theorem.} {\it A Banach space is extremely projective in {\bf Ban} if and only if it is isometrically isomorphic to $l_1(\Lambda)$ for some index set $\Lambda$.}

\medskip
(This theorem has a substantial operator space ( = `non-commutative') version; see~\cite[Thm. 3.14]{ble})

\medskip
With obvious modifications, one can define spaces, {\it topologically projective in {\bf Nor}} and
{\it topologically projective in {\bf Ban}}.

\medskip
The discussion of extremely projective spaces in {\bf Nor} and in {\bf Ban}, including the proof of the formulated proposition, will be presented in Section 2. In particular, we shall see that the same space can be extremely projective in {\bf Ban} but not extremely, and even topologically, projective in {\bf Nor} (Proposition 2.5).

\bigskip
From just spaces let us turn again to modules.

\bigskip
The interest in extreme versions of basic homological notions was stimulated by the fundamental 
Arveson--Wittstock Theorem in quantum functional analysis, playing in that area the role of Hahn--Banach 
Theorem in classical functional analysis. First, extremely flat and extremely injective modules appeared 
(in~\cite{he4,wit}), however only for the special case of some highly non--commutative operator algebras 
and for the class of the so--called semi--Ruan modules in the capacity of the distinguished category $\f$. 
The results of cited papers led to some extension theorems of Hahn--Banach type for modules. These theorems, 
in their turn, led to a conceptually new proof of the Arveson--Wittstock Extension Theorem in its 
non--coordinate presentation. This was done in~\cite{he4} for the initial case of that theorem, dealing 
with operator spaces, and in~\cite{wit} for the more sophisticated case of operator modules.

In the present paper we concentrate on projectivity and consider another class of  `popular' algebras, which is opposite, in a sense, to those in~\cite{he4,wit}. We mean commutative normed algebras, in fact algebras, consisting of sequences. These algebras apparently represent the next degree of complication after $\C$. Nevertheless we hope to show that even in this case, after the proper choice of $\f$, there is something to be said.

 Denote by ${\bf p}^n$ the sequence $(0,\dots,0,1,0,0\dots)$ with 1 as its $n$-th term,
and by $c_{00}$ the linear space of finite sequences, that is $span\{{\bf p}^n; n=1,2,\dots\}$.

\medskip
{\bf Definition II.} Let $A$ be a normed algebra, consisting of some complex--valued sequences and equipped with the coordinate--wise operations. We say that $A$ is a {\it sequence algebra}, if it
contains, as a dense subalgebra, $c_{00}$, and for all $n$ we have $\|{\bf p}^n\|=1$.

\medskip
We see that the class of sequence algebras includes $c_0$, all $l_p; 1\le p<\ii$ (but not $l_\ii$),
Fourier algebras of discrete countable groups (after rearranging the respective domains as sequences), and many other algebras.

\bigskip
The main results of this paper are Theorem I and II below. They give, within a certain substantial class of normed modules over a sequence algebra, a full characterization of extremely projective modules with respect to that class. Theorem I essentially deals with, generally speaking, non--complete modules. (The interest to the `non--Banach' case is derived from the cited papers on the extreme flatness). However, from this theorem its Banach counterpart
easily follows. This is Theorem II, which in one important point sounds  essentially different.

We proceed to define, after some preparatory observations, the distinguished class of modules, playing the role of $\f$ (see above).

Let $A$ be a sequence algebra, $X$ a normed $A$-module. Often, when there is no danger of confusion, for $x\in X$ we shall write $x_n$ instead of ${\bf p}^n\cd x$ and call the latter the {\it $n$--th coordinate} of $x$. Of course, we have ${\bf p}^n\cd x_n=x_n$. Further, we set $X_n:=\{{\bf p}^n\cd x; x\in X\}$ for every $n\in\N$. We see that $X_n$ is a subspace of $X$. Moreover, $X_n$ is obviously a submodule with the outer multiplication
$$
a\cd x=a_nx; a=(a_1,\dots,a_n,\dots)\in A, x\in X_n.\eqno(1)
$$
It will be called the {\it $n$--th coordinate subspace (or submodule)} of $X$.

\medskip
{\bf Definition III}. An $A$-module $X$ is called {\it homogeneous} if, for every $x,y\in X$, the inequalities $\|x_n\|\le\|y_n\|; n\in\N$ imply that $\|x\|\le\|y\|$.

\medskip

It immediately follows that the equalities $\|x_n\|=\|y_n\|; n=1,2,\dots$ for any elements $x,y$ in a homogeneous module imply that $\|x\|=\|y\|$. Thus in a homogeneous module the norm of an element is completely determined by the norms of its coordinates.


\medskip
For many typical sequence algebras the class of homogeneous modules is fairly wide. In particular, it is easy to show that all essential normed modules over $c_0$, consisting of complex--valued sequences, are homogeneous.
 Besides, $l_p$-sums; $1\le p\le\ii$ of arbitrary families of normed spaces are obviously homogeneous $A$-modules.
 (In both examples we mean the coordinate--wise outer multiplication).

On the other hand, it is obvious that a homogeneous normed $A$-module $X$ is necessarily faithful, 
that is for every $x\in X$ the equality $a\cd x=0$ for all $a\in A$ implies $x=0$.

In the following two theorems $A$ is an arbitrary sequence algebra. We denote by ${\cal H}$ the full 
subcategory in $A-mod$, consisting of all essential homogeneous modules, and by $\overline{{\cal H}}$ the full subcategory 
in ${\cal H}$, consisting of Banach modules.

We call an element $x$ of a given normed $A$--module $X$ {\it finite}, if the sequence of its
coordinates is finite, that is $x_n=0$ for all sufficiently large $n$. We denote by $X_{00}$ the submodule of $X$, consisting of finite elements. If $X$ is faithful, in particular, homogeneous, then  for  every $x\in X_{00}$
  we obviously have
$$
x=\sum_nx_n,
$$
(needless to say that the latter sum is well defined), and $X_{00}$ is exactly the {\it algebraic} direct sum of coordinate submodules $X_n; n\in\N$.

We say that $X$ {\it is of finite type}, or $X$ {\it has finite type}, if $X=X_{00}$.


\medskip
{\bf Theorem I}. {\it A module $X\in{\cal H}$ is extremely projective with respect to ${\cal H}$ if and only if it satisfies the following two conditions:

(i) For every $n\in\N$, the $n$--th coordinate subspace $X_n$ is
extremely projective in {\bf Nor}, or equivalently} (see Proposition above), {\it $X_n$ is a near--retract of $l_1^0(\Lambda_n)$ for some index set $\Lambda_n$.

(ii) $X$ is of finite type.

Moreover, if $X$ is at least topologically projective, (ii) is again valid. }

\medskip
The Banach counterpart of the formulated theorem is

\medskip
{\bf Theorem II}. {\it A module $X\in\ov{{\cal H}}$ is extremely projective with respect to $\ov{{\cal H}}$ if and only if for every $n\in\N$ the $n$--th coordinate subspace $X_n$ is extremely projective as a Banach space, or equivalently} (see above), {\it $X_n$ is isometrically isomorphic to $l_1(\Lambda_n)$  for some index set $\Lambda_n$.}

\medskip
Thus, speaking informally, in both theorems the answer depends not on the norm on the
whole module but only on the norms of its coordinate subspaces. In particular, all Banach
homogeneous modules, consisting of sequences,  are extremely projective with respect to $\ov{{\cal H}}$, 
but neither of them is extremely projective with respect to ${{\cal H}}$, when it is infinite--dimensional. 
On the other hand, submodules of these modules, 
consisting of finite sequences, are extremely projective with respect to ${{\cal H}}$.


\bigskip
\centerline{{\bf 1. Preparatory observations}}

\bigskip
First, we consider the general case of an arbitrary normed algebra $A$ and an arbitrary distinguished
full subcategory $\f$ in $A-mod$.

\medskip
{\bf Proposition 1.1.} {\it Let $Q$ be a normed $A$-module, and  $P$ its near--retract. Assume that $Q$ is extremely projective with respect to $\f$. Then the same is true for $P$.}

\smallskip
PROOF. Suppose we are given $\tau$, $\va$ and $\e$ as in Definition I. Fix any $\de>0$ such that $\|\va\|\de+\de+\de^2<\e$. By the assumption, for $\va_0:=\va\s$  there exists its lifting $\psi_0:Q\to Y$ such 
that $\|\psi_0\|<\|\va_0\|+\de$, and there exist bounded morphisms $\s:Q\to P$ and $\rho:P\to Q$ such that $\s\rho=\id_P$ and $\|\rho\|<1+\de$. Now we see that we $\psi:=\psi_0\rho$ is a lifting we need. $\Box$

\bigskip
If a distinguished category $\f$ is given, we denote by $\overline{\f}$ the full subcategory of $A-mod$ whose objects are Banach $A$-modules that are completions of modules from $\f$.

\medskip
{\bf Proposition 1.2.} {\it Suppose that for every module in $\f$ its completion also belongs to $\f$.
Then, for every $P\in A-mod$, which  is extremely projective with respect to $\f$, its completion
$\ov{P}$ is extremely projective with respect to $\ov\f$.}

\smallskip
PROOF. Suppose we are given a coisometric morphism $\tau:Y\to X; X,Y\in\ov{\f}$, a bounded morphism  $\va:\ov{P}\to X$ and $\e>0$. Consider the restriction $\va_0$ of $\va$ to $P$. By the assumptions on $\f$ and $P$, there exists a lifting $\psi_0:P\to Y$ of $\va_0$ with $\|\psi_0\|<\|\va_0\|+\e$. Since $Y$ is a Banach module, $\psi_0$ has the continuous extension $\psi:\ov{P}\to Y$, which is obviously a lifting of $\va$. Since $\|\va_0\|=\|\va\|$ and $\|\psi\|=\|\psi_0\|$, the previous estimate of $\|\psi_0\|$ gives the desired estimate of $\|\psi\|$. $\Box$

\bigskip
In the remaining part of this section  we concentrate on the case when $A$ is a sequence algebra. This material will be used in Section 3.

Note that an arbitrary morphism $\va:Y\to X$ in $A-mod$ gives rise, in an obvious way, to the sequence of its birestrictions between the respective coordinate submodules (cf. Introduction). These morphisms will be denoted by $\va_n:Y_n\to X_n; n=1,2,\dots$ and called {\it coordinates submorphisms of $\va$}.

\medskip
{\bf Proposition 1.3.} {\it If $\tau:Y\to X$ is a coisometric morphism in $A-mod$, then, for every $n$, $\tau_n:Y_n\to X_n$ is also a coisometry.}

\smallskip
PROOF. Take $x\in X_n$ and $\e>0$. By assumption, there exists $\widetilde y\in Y$ with $\tau(\widetilde y)=x$ and
 $\|\widetilde y\|<\|x\|+\e$. Then for $y:={\bf p}^n\cd\widetilde y$ we obviously have $\tau(y)=x$ and $\|y\|\le\|\widetilde y\|$.  The rest is clear. $\Box$

\medskip
Now fix, for a moment, some $n\in\N$ and denote by $\f_n$ the full subcategory in {\bf Nor}, whose objects are $n$-th coordinate subspaces of modules in our distinguished category $\f$.

\medskip
{\bf Proposition 1.4.} {\it Let $\f$ contains, with its every module, the $n$-th coordinate submodule of the latter. Assume that a normed $A$-module $P$ is extremely projective with respect to $\f$. Then $P_n$ is extremely projective in {\bf Nor} with respect to $\f_n$.}

\smallskip
PROOF. Take an arbitrary $F,E\in\f_n$, a coisometric operator $\tau:F\to E$, a bounded operator $\va:P_n\to E$ and $\e>0$. We must find a lifting operator $\psi:P_n\to F$ such that $\|\psi\|<\|\va\|+\e$.

Our $E$ and $F$ are underlying spaces of $n$-th coordinate submodules, say $X$ and $Y$, of some modules in $\f$.
By the assumption,
$X,Y\in\f$, and (1)
implies that $\tau$, as a map between these modules, is a module morphism.

Consider the map $\widetilde\va:P\to X:x\mt\va(x_n)$. It follows from the same (1) that it is a morphism of $A$-modules. Besides, we obviously have $\|\widetilde\va\|=\|\va\|$. Therefore, the assumption on $P$
provides a lifting $\widetilde\psi:P\to Y$ of $\widetilde\va$ such that $\|\widetilde\psi\|<\|\va\|+\e$.

Denote by $\psi$ the restriction of $\widetilde\psi$ to $P_n$. Then we easily see that
it is a lifting of $\va$, and $\|\psi\|\le\|\widetilde\psi\|$. The rest is clear. $\Box$

\medskip
Let us turn to the special properties of homogeneous modules.

\medskip
 For every $N=1,2,\dots$ we set ${\bf P}^N:=\sum_{n=1}^N{\bf p}^n\in A$.

\medskip
{\bf Proposition 1.5}. {\it If $A$-module $X$ is essential and homogeneous, then for every  $x\in X$  we have}
$$
x=\lim_{N\to\ii}{\bf P}^N\cd x.
$$

\smallskip
PROOF. Take $x$ and $\e>0$. It easily follows from Definition II that there is $y\in Y$ of the form
$\sum_{k=1}^na^k\cd z^k; a^k\in c_{00}, z^k\in X$ such that $\|x-y\|<\e/2$. For all $N\in\N$ we have
$$
\|x-{\bf P}^N\cd x\|\le\|x-y\|+\|y-{\bf P}^N\cd y\|+\|{\bf P}^N\cd y-{\bf P}^N\cd x\|.
$$
But, because of the choice of $y$, for some $M\in\N$ we have $y={\bf P}^N\cd y$
for all $N>M$. Besides, the homogeneity of $X$ implies that $\|{\bf P}^N\cd(y-x)\|\le\|y-x\|$. Therefore for all
$N>M$ we have $\|x-{\bf P}^N\cd x\|<\e$. $\Box$

%
%

\medskip
{\bf Proposition 1.6}. {\it If  $A$-module $X$ is essential and homogeneous, the same is true for its completion $\overline{X}$.} (In other words, $X\in{\cal H}$ implies $\ov{X}\in\ov{{\cal H}}$)

\smallskip
PROOF. It is immediate that $\ov{X}$ is essential. Let us prove that it is homogeneous. Take $x,y\in\overline{X}$ with $\|x_n\|\le\|y_n\|$ for all $n$. We must show that $\|x\|\le\|y\|$. By virtue of the previous proposition
we can assume, without loss of generality, that for some natural $N$ we have $x_n=y_n=0$ for all $n>N$.

Choose sequences $x^k,y^k\in X;k=1,2,...$ converging to $x$ and $y$, respectively.
Using Proposition 1.5, we can assume,
again without loss of generality, that $x_n^k=0$ whenever $x_n=0$.

For all $n$ we have $\lim_{k\to\ii}\|x^k_n\|=\|x_n\|$ and $\lim_{k\to\ii}\|y^k_n\|=\|y_n\|$.
Taking into account what we have just assumed, we see that
for every $\e>0$ there exists a natural $M$ such that for all $k>M$ and for all $n$ we have $\|x^k_n\|\le(1+\e)\|y^k_n\|$.

But $X$ is homogeneous, and therefore for the same $k$ we have $\|y^k\|\le\|(1+\e)x^k\|$.
It remains to pass to limits. $\Box$

\medskip
Now we recall a property of Banach (but not just normed) spaces, which is very well known and which is
easy to prove. Namely, {\it if $\sum_{n=1}^\ii x_n$ is a converging series in a Banach space $E$,
then there exists another series in $E$, say $\sum_{n=1}^\ii\bar x_n$, such that, for some sequence
$\lm_n\in{\B R}; \lm_n\ge1$ and $\lm_n\to\ii$, we have $\bar x_n=\lm_nx_n$ and this new series is
still converging.} From this one can easily deduce

\medskip
{\bf Proposition 1.7}. {\it Let $X$ be essential and homogeneous. Then for every $x$ in its completion $\overline{X}$ (in particular, in $X$) there exists $\bar x\in\overline{X}$ such that $\bar x_n=\lm_n^xx_n$, where $\lm_n^x\in{\B R}; \lm_n^x\ge1$ and $\lm_n^x\to\ii$.}

\medskip
Speaking informally, the coordinates of $\bar x$, being proportional to those of $x$, tend to 0 essentially slower.

\smallskip
PROOF. By the assumptions on $X$, the limit equality in Proposition 1.5 can be rewritten as $x=\sum_{n=1}^\ii x_n$, the sum of the series of its coordinates. Take $\bar x_n,\lm_n$ as before and rewrite $\lm_n$ as $\lm_n^x$. Then
$\bar x:=\sum_{n=1}^\ii\bar x_n$ obviously fits. $\Box$

\bigskip
At the end of the section we describe a certain way to construct, with the help of a given homogeneous
 module of finite type, some new homogeneous modules.

Suppose that we have an algebraic $A$--module $X$ of finite type. Set $M_X:=\{n\in\N:X_n\ne0\}$ and denote by $c_{00}^+(X)$ the set (actually, a cone) of all finite non-negative sequences $\xi=(...,\xi_n,...)$ such that $\xi_n=0$ whenever $n\notin M_X$.

Now assume that our $X$ is a normed homogeneous $A$--module. Introduce the function $f_X:c_{00}^+(X)\to\R$,
taking $\xi$ to $\|x\|$, where $x$ is an (obviously existing) element of $X$ such that $\|x_n\|=\xi_n$.
By homogeneity of $X$, this function is well defined.

\medskip
{\bf Proposition 1.8.} {\it The function $f_X$ has the following properties:

(i)  If $\xi\in c_{00}^+(X)$ is not zero, then $f_X(\xi)>0$.

(ii) If $\xi\in c_{00}^+(X)$ and $\lm>0$, then $f_X(\lm\xi)=\lm f_X(\xi)$

(iii) If $\xi,\eta\in  c_{00}^+(X)$, and $\xi\le\eta$, then $f_X(\xi)\le f_X(\eta)$

(iv) If $\xi\in c_{00}^+(X)$ and $a\in A$, then $f_X(|a|\xi)\le\|a\|f_X(\xi)$.

(v) If $\xi,\eta\in c_{00}^+(X)$, then $f_X(\xi+\eta)\le f_X(\xi)+f_X(\eta)$

(vi) If $n\in M_X$, then $f_X({\bf p}^n)=1$.}

\smallskip
PROOF. The properties (i)-(iv) are immediate.

\smallskip
(v) Take $x\in X$ such that $\|x_n\|=(\xi+\eta)_n$ for all $n$ and thus $f_X(\xi+\eta)=\|x\|$.
If we have, for a given $n$, $\xi_n+\eta_n>0$, then we set $\lm_n:={\xi_n}/({\xi_n+\eta_n})$ and $\mu_n:={\eta_n}/({x_n+\eta_n})$;
otherwise we set $\lm_n=\mu_n=0$. After this, we set $y:=\sum_n\lm_nx_n$ and $z:=\sum_n\mu_nx_n$. Then, of course, 
we have $f_X(\xi)=\|y\|$ and $f_X(\eta)=\|z\|$. But $x=\sum_nx_n=y+z$ and hence $\|x\|\le\|y\|+\|z\|$. The desired inequality follows.

\smallskip
(vi) Take $x\in X_n$ of norm 1. Then the sequence $(\|x_1\|,\|x_2\|,...)$ is exactly ${\bf p}^n$. The rest is clear. $\Box$

\medskip
Sometimes we shall refer to $f_X$ as to the {\it function, associated with the module $X$.}

\medskip
We began with a module and got a function. Now we proceed in the opposite direction.

Let $X$ be an algebraic $A$--module $X$ of finite type such that, for $x\in X$, the equalities ${\bf p}^n\cd x=0; n\in\N$ imply that $x=0$. Suppose that
every coordinate subspace $X_n$ is equipped with a norm, say $\|\cd\|_n$. Fix an arbitrary
function  $f:c_{00}^+(X)\to\R$, possessing the properties (i)-(v) of the previous proposition. Now for $x\in X$ we set $\|x\|:=f(\xi)$, where $\xi_n:=\|x_n\|_n$.

\medskip
{\bf Proposition 1.9.} {\it The assignment $x\mt\|x\|$ is a norm on $X$, making it a normed homogeneous $A$--module. If, in addition, our function has the property (vi), then, for every $n$, the new norm $\|\cd\|$, being restricted on $X_n$, coincides with the initial norm $\|\cd\|_n$.}

\smallskip
PROOF. Of the standard properties of a norm, only the triangle inequality is not immediate. Take $x,y,z\in X$ such that $x=y+z$.
Take the sequences $\xi,\eta,\zeta\in c_{00}^+(X)$ such that $\xi_n:=\|x_n\|_n, \eta_n:=\|y_n\|_n$ and $\zeta_n:=\|z_n\|_n$. Then we obviously have $\xi\le\eta+\zeta$, and the properties (iii) and (v) of $f$ imply that
$f(\xi)\le f(\eta+\zeta)\le f(\eta)+f(\zeta)$. Consequently, $\|x\|\le\|y\|+\|z\|$.

Now take $a\in A$ and $x\in X$. Then, for the sequences $\xi$ and $\eta$ with $\xi_n:=\|x_n\|_n$ and
$\eta_n:=\|(a\cd x)_n\|_n$, we have $\eta=|a|\xi$. Therefore the property (iv) of $f$ implies that $\|a\cd x\|\le\|a\|\|x\|$.

 Thus $X$ became a normed module, which is obviously homogeneous.

 Finally, note that if, for some $n$, we have $x\in X_n$, then the sequence \\ $(\|x_1\|_1,\|x_2\|_2,\dots)$ coincides with $\|x_n\|_n{\bf p}^n$, and therefore $\|x\|=\|x_n\|_nf({\bf p}^n)$. The last assertion follows. $\Box$

 \medskip
If the initial algebraic module $X$
with the indicated properties together with norms $\|\cd\|_n;n\in\N$ are fixed, we shall denote  the constructed homogeneous module by $X^f$.


 \medskip
{\bf Proposition 1.10.} {\it (i) If $X$ is an arbitrary homogeneous module of finite type, then $X^{f_X}=X$

 (ii) If we have the data of Proposition 1.9, and $f$ satisfies the properties (i)-(vi) of Proposition 1.8, then $f_{X^f}=f$.}

\smallskip
PROOF. (i) is immediate, and (ii) follows from the last assertion in Proposition 1.9. $\Box$

 \medskip
Finally, suppose that we have two $A$--modules in the pure algebraic sense, say $X$ and $Y$, and, for every $n$, a linear operator $\va_n:X_n\to Y_n$ is given. Suppose, further, that $X$ is of finite type. Then there exists, for every $x\in X$, the well defined element $\va(x):=\sum_n\va_n(x_n)\in Y$. In this way we obtain a map $\va:X\to Y$, which is, of course, a morphism. We shall call it {\it morphism, generated by operators $\va_n$.} Note that for every $x$ and $n$ we obviously have
$$
\va(x)_n=\va_n(x_n).\eqno(2)
$$

\medskip
{\bf Proposition 1.11.} {\it Let $X$ and $Y$ be normed homogeneous modules of finite type with $M_X=M_Y$ (that is, $c_{00}^+(X)=c_{00}^+(Y)$) and with the same associated function.
Suppose that for every $n$ we are given a bounded operator $\va_n:X_n\to Y_n$, and $C:=\sup\{\|\va_n\|; n=1,2,\dots\}<\ii$. Then the morphism $\va$, generated by operators $\va_n$, is bounded, and $\|\va\|=C$.}

\smallskip
PROOF. Take an arbitrary $x\in X$ and denote by $\xi$ and $\eta$ the sequences with $\xi_n:=\|x_n\|$ and $\eta_n:=\|\va(x)_n\|$. It follows from (2) and the assumption on $\va_n$ that $\eta_n\le C\xi_n$. Since
$f_X=f_Y$, Proposition 1.8(iii,ii) implies $f_Y(\eta)\le Cf_X(\xi)$, that is, by definition of associated functions,  $\|\va(x)\|\le C\|x\|$. Thus we have $\|\va\|\le C$. Further, we see from (2) that, for every $n$, the restriction of $\va$ to $X_n$ is $\va_n$. From this we have $C\le\|\va\|$. $\Box$


\bigskip
\centerline{{\bf 2. On extremely projective non-complete and complete normed spaces}}

\bigskip
Looking for extremely projective spaces in the `non--complete context', one inevitably pays attention to the
spaces $l_1^0(\Lambda)$, mentioned in the Introduction. For a moment, fix a set $\Lambda$ and for every $\nu\in\Lambda$ denote by $e_\nu$ the function in $l_1^0(\Lambda)$, taking $\nu$ to 1 and other
points of $\Lambda$ to 0. Clearly, the set $e_\nu;\nu\in\Lambda$ is a linear basis in $l_1^0(\Lambda)$. In what
follows, such a basis will be called {\it natural}.

The following two statements must be well known, at least as a math folklore.

\medskip
{\bf Proposition 2.1.} {\it Let $E$ be a normed space, $\al:\Lambda\to E$ a map with a bounded image. Then there exists a bounded operator $\va:l_1^0(\Lambda)\to E$, uniquely defined by $e_\nu\mt\al(\nu); \nu\in\Lambda$. Moreover, $\|\va\|=\sup\{\|\al(\nu)\|;\nu\in\Lambda\}$.}

\smallskip
PROOF. The assertion obviously follows from the nature of the norm on $l_1^0(\Lambda)$. $\Box$

\medskip
{\bf Proposition 2.2.} {\it The space $l_1^0(\Lambda)$, where $\Lambda$ is an arbitrary index set, is extremely projective in {\bf Nor}.}

\smallskip
PROOF. Take a coisometry $\tau:F\to E$ in {\bf Nor}, a bounded operator $\va:l_1^0(\Lambda)\to E$ and $\e>0$. Fix, for every $\nu\in\Lambda$, an arbitrary vector $y_\nu\in F$ such that $\tau(y_\nu)=\va(e_\nu)$ and $\|\tau(y_\nu)\|<\|\va(e_\nu)\|+\e$. Consider the operator $\psi:l_1 ^0(\Lambda)\to F$, well defined by $e_\nu\mt y_\nu$. Evidently, it is a lifting of $\va$ with respect to $\tau$. Besides, the previous proposition implies that $\|\psi\|<\|\va\|+\e$. $\Box$

\medskip
Let $E$ be a normed space. Denote its unit sphere by $S_E$. Consider the normed space $l_1^0(S_E)$ (i.e. $l_1^0(\Lambda)$ with $S_E$ in the role of the relevant index set) and its natural linear basis $e_x; x\in S_E$ (see above). Then Proposition 2.1 provides a bounded operator $\tau_E:l_1^0(S_E)\to E: e_x\mt x; x\in S_E$; obviously, it is a coisometry. We shall call it the {\it canonical coisometry for $E$}.

\medskip
Now we can prove the proposition, formulated in the Introduction.

\medskip
PROOF. The `if' part follows from the previous proposition and Proposition 1.1; in the latter we set $A-mod:=\f=:{\bf Nor}$ and $Q:=l_1^0(\Lambda)$. Conversely, suppose that a certain $P$ is extremely projective in {\bf Nor}. Evidently, we can assume that $P\ne0$. Consider the canonical coisometry $\tau_P:l_1^0(S_P)\to P$ for $P$.
Then, by the assumption on $P$, for every $\e>0$ there exists a lifting $\rho:P\to l_1^0(B)$ of the identity operator $\id_P$ on $P$ such that $\|\rho\|<\|\id_P\|+\e$. The rest is clear. $\Box$

\medskip
Of course, a much more sound and transparent statement would be:

  a normed space is extremely projective as a normed space if and only if it is isometrically isomorphic to $l_1^0(\Lambda)$ for some index set $\Lambda$.

  But is it true? We do not know.

The question seems reasonable, especially because it has, as a background, the Grothendieck Theorem that was formulated in the Introduction.

\medskip
{\small {\bf Remark.} In literature, speaking about the latter theorem, they usually cite~\cite{gro} (see, e.g.,~\cite[p. 182]{heb}). Basically, it is correct. At the same time, despite the paper~\cite{gro} contains all needed ingredients for the proof, the theorem itself is not explicitly formulated
 (could it be that the author just thought it unnecessary ?)
 By `ingredients' we mean the following two statements, formulated (needless to say, in equivalent terms) and completely proved:

\medskip
(i) A Banach space, which is topologically projective in {\bf Ban} and isometrically isomorphic to some $L_1(\Omega,\mu)$, is isometrically isomorphic to $l_1(\Lambda)$~\cite[Prop. 2]{gro}.

(ii) A Banach space $F$ is isometrically isomorphic to some $L_1(\Omega,\mu)$
 if and only if, for every coisometry, {\it which is also an adjoint operator}, say $\tau^*:Y^*\to X^*$, and every bounded operator $\va:F\to X^*$, there is a lifting $\psi$ of $\va$ with respect to $\tau^*$ such that $\|\psi\|=\|\va\|$ \\ \cite[Prop. 1(3)]{gro}. (Note the absence of any $\e$ -- A.H.)

\medskip
After these statements, very little remains to complete the proof. Indeed, take our given projective space, say $F$, and fix for a time $\tau^*$ and $\va$ as above.
Then for every $\e>0$ there is a lifting $\psi_\e:E\to Y^*$ of $\va$ with $\|\psi_\e\|<\|\va\|+\e$. But the range of all these $\psi_\e$ is a dual space. Therefore, using the standard argument, based on the Banach-Alaoglu Theorem, we can take a cluster point, say $\psi$, of the set $\psi_\e; \e>0$ in the weak$^*$ topology of $\bb(E,Y^*)=(E\widehat\otimes Y)^*$. It is easy to see that $\psi$ has the property,
indicated in (ii). Thus $F$ satisfies (i), and we are done.}

\medskip
We return to general normed spaces. Combining the Grothendieck Theorem with Proposition 1.2, considered for the particular case $A-mod=\f={\bf Nor}$, we immediately get

\medskip
{\bf Proposition 2.3.} {\it Suppose $P$ is extremely projective in {\bf Nor}. Then it is, up to an isometric isomorphism, a dense subspace in $l_1(\Lambda)$ for some index set $\Lambda$.}

\medskip
Thus the question, formulated above, can be posed in the following, somewhat more detailed form. For a given $\Lambda$, which dense subspaces of $l_1(\Lambda)$ are extremely projective in {\bf Nor} (like $l_1^0(\Lambda)$), and which are not?

Concluding this section, we point out another necessary condition of the property under discussion.

\medskip
{\bf Proposition 2.4.} {\it Let $P$ be a separable normed space. Suppose that it is extremely, or at least topologically projective in {\bf Nor}. Then it has at most countable linear dimension.}

\smallskip
PROOF.
Consider the canonical coisometry  $\tau_P:l_1^0(S_P)\to P$ for $P$. Then, by the second assumption,
the identity operator $\id_P$ has a bounded lifting $\psi:P\to l_1^0(S_P)$ with respect to $\tau$.

Now, using the first assumption, choose an arbitrary dense subset, say $\{x_n;  \\ n\in N\}$, in $P$. Since $l_1^0(S_P)$ has the {\it linear} basis $e_x; x\in S_P$, every $\psi(x_n); n\in\N$ has a form $\sum_{k=1}^{m_n}\lm^{(n)}_ke_{x_k}^n$ for some ${m_n}\in\N, x_k^n\in S_P$ and $\lm^{(n)}_k\in\C; k=1,\dots,m_n$. Therefore all $\psi(x_n); n\in\N$ belong to the linear span of of all $e_{x_k^n};n\in\N,k=1,\dots,m_n$, denoted for brevity by $F$. The latter is obviously closed in $l_1^0(S_P)$, the set $\{x_n; n\in\N\}$ is dense in $P$, and $\psi$ is continuous. Consequently, $\psi$ maps $P$ into $F$. Further, $\psi$, being a right inverse to some map (namely, $\tau$), is injective. Therefore, $\psi$ implements a linear isomorphism of $P$ onto its image in $F$, and this image, of course, has at most countable linear dimension. $\Box$

\medskip
From this, one can easily deduce

\medskip
{\bf Proposition 2.5.} {\it Let $\Lambda$ be an arbitrary infinite index set. Then the space $l_1(\Lambda)$}
(being, thanks to Grothendieck, extremely projective in {\bf Ban}) {\it is not topologically, and hence extremely, projective in {\bf Nor}.}

\smallskip
PROOF. It is obvious that $l_1(\Lambda)$ has $l_1(\N)$, that is $l_1$, in the capacity of its near--retract (actually, of its retract in the category of normed spaces and contracting operators). Therefore, by virtue of Proposition 1.1, it is sufficient to prove our assertion for the latter space. But $l_1$ is, of course, separable, and its linear dimension, by the old theorem of L\"ovig (see, e.g.,\cite{loe}) is continuum. Now the previous proposition works. $\Box$

\bigskip
\centerline{{\bf 3. The proof of Theorems I and II}}

\bigskip
Up to the end of the paper, $A$ is a fixed sequence algebra, the word `module' always  means `$A$--module`, and the word `morphism' always  means `morphism of $A$--modules'.

The following theorem is the main step in the proof of the `if` part of \\ Theorem I.

\medskip
{\bf Theorem 3.1.} {\it Let $P$ be a homogeneous normed $A$-module of finite type such that, for every $n,\; P_n=l_1^0(\Lambda_n)$ for an index set $\Lambda_n$. Then $P$ is extremely projective with respect to ${\cal H}$.}

\smallskip
PROOF. For every $n$ we denote by $e^n_\nu;\in\Lambda_n$ the natural basis in $l_1^0(\Lambda_n)$, by $S$ the unit sphere in $P$ and by $E(S)$ the set of such $z$ in $S$ that, for every $n,\; z_n$ is a multiple of $e^n_\nu$ for some $\nu\in\Lambda_n$.

\medskip
We need two preparatory assertions. In both of them $X$ is an arbitrary normed module, and $\va:P\to X$ is a bounded morphism.

\medskip
{\bf Lemma.} {\it For every $x\in P$ and $N\in\N$ there exists $x^{(N)}\in P$ such that $x^{(N)}_N$ is a multiple of some $e^N_\nu;\nu\in\Lambda_N$, $x^{(N)}_n=x_n$ for all $n\ne N$, $\|x^{(N)}\|=\|x\|$, and $\|\va(x^{(N)})\|\ge\|\va(x)\|$.}

\medskip
Proof of the lemma. Of course, we can assume that $x_N\ne0$, and thus it has a form $\sum_{k=1}^m\lm_ke^{N}_{\nu_k}$, where $\lm_k$ are non-zero complex numbers, $\nu_k\in\Lambda_N$. Set, for brevity, $\mu_k:=|\lm_k|/\|x_N\|$ and consider, for every $k=1,\dots,m,$ the element $x^k\in P$ with the coordinates $x^k_N:=\mu^{-1}_k\lm_ke^{N}_{\nu_k}$ and $x^k_n:=x_n$ for $n\ne N$. Obviously, for all $k$  we have $\|x^k_N\|=\|x_N\|$. Since $P$ is homogeneous, this implies $\|x^k\|=\|x\|$. Further, since $\sum_{k=1}^m|\lm_k|=\|x_N\|$, we have $\sum_{k=1}^m\mu_kx^k=x$. Therefore $\sum_{k=1}^m\mu_k\va(x^k)=\va(x)$, and we see that $\va(x)$ is a convex combination
of elements $\va(x^k)$ in $X$. Consequently, for at least one of $k$ we have $\|\va(x^k)\|\ge\|\va(x)\|$.
Thus the respective $x^k$ has all properties of the desired $x^{(N)}$.

\medskip
{\bf Lemma.} {\it We have $\|\va\|=\sup\{\|\va(z)\|; z\in E(S)\}$.}

\smallskip
Proof of the lemma. Take an arbitrary $x\in S$. Apply to this elements the previous proposition for the case $N:=1$. Then apply to the resulting element $x^1$ (also belonging, by the construction, to $S$) the same proposition, this time for the case $N:=2$, and so on. Since $x$ is finite, eventually we come to an element, say $z$, belonging to $E(S)$ and such that $\|\va(z)\|\ge\|\va(x)\|$. The rest is clear.

\newpage

\medskip
{\bf The end of the proof of Theorem 3.1.}

\smallskip
Suppose we are given modules $X,Y\in{\cal H}$, a coisometric morphism $\tau:Y\to X$, a bounded morphism $\va:P\to X$, and $\e>0$. Consider the coordinate submorphisms $\va_n:P_n\to X_n, \tau_n: Y_n\to X_n$ and choose $\de>0$ with $(1+\de)^2\|\va\|<\|\va\|+\e/2$. By Proposition 1.3, all $\tau_n$ are coisometries. Therefore for every $\nu\in\Lambda_n$ there exists an element $y^n_\nu\in Y_n$ such that $\tau_n(y^n_\nu)=\va(e^n_\nu)$ and $\|y^n_\nu\|\le\|\va(e^n_\nu)\|+\de\|\va(e^n_\nu)\|$.
Besides, Proposition 2.1 provides, for every $n$, an operator $\psi_n:P_n\to Y_n$, well defined by $e^n_\nu\mt y_{n_\nu}; \nu\in\Lambda_n$. Consider the morphism $\psi:P\to Y$, generated (see above) by the operators $\psi_n$. Clearly, it is a lifting of $\va$ with
respect to $\tau$.

Now take $x\in E(S)$. For every $n$ we have $x_n=\lm^n_\nu e^n_\nu$ for some $\nu\in\Lambda_n$ and $\lm^n_\nu\in\C$. Hence $\va_n(x_n)=\lm^n_\nu\va(e^n_\nu)$ and $\psi_n(x_n)=\lm^n_\nu y^n_\nu$. Therefore we have

$$
\|\psi_n(x_n)\|\le(1+\de)\|\va_n(x_n)\|.\eqno(3)
$$

But, since $\tau$ is a coisometry, there exists $z\in Y$ with $\tau(z)=\va(x)$ and $\|z\|\le\|\va(x)\|+\de\|\va(x)\|$.
Besides, we have $\tau(z_n)=\tau(z)_n=\va(x)_n=\va_n(x_n)$ and hence $\|\va_n(x_n)\|\le\|z_n\|$. Together with (3), this gives $\|\psi(x)_n\|=\|\psi_n(x_n)\|\le(1+\de)\|z_n\|$. From this, by homogeneity, we obtain that
$$
\|\psi(x)\|\le(1+\de)\|z\|\le(1+\de)^2\|\va(x)\|\le(1+\de)^2\|\va\|<\|\va\|+\e/2.
$$
 Therefore the previous lemma gives $\|\psi\|\le\|\va\|+\e/2$. The rest is clear. $\Box$

\bigskip
{{\bf The proof of the `if' part of Theorem I}}.

\bigskip
Now the coordinate subspaces $P_n$ of our homogeneous module $P$ of finite type are arbitrary spaces that are extremely projective in {\bf Nor}. By virtue of the proposition, formulated in the Introduction and proved in Section 2, for every $n\in\N$ there exists an index set $\Lambda_n$ and an operator $\s_n:l_1^0(\Lambda_n)\to P_n$ which is a
near-retraction.

Set, for brevity, $Q_n:=l_1^0(\Lambda_n)$
 (thus, in particular, $Q_n=0$ exactly when $P_n=0$). Denote by $Q$ the algebraic sum $\bigoplus_nQ_n$. Clearly, $Q$ is an $A$-module of finite type with respect to the
outer multiplication, well defined, for $a\in A$ and $x\in Q$ of the form $x=\sum_nx_n; x_n\in Q_n$, by $a\cd x:=\sum_na_nx_n$.

\smallskip
Of course, $c_{00}^+(Q)=c_{00}^+(P)$, and hence we can endow $Q$ with a norm by the recipe of Proposition 1.9, taking, in the capacity of $f$, the function $f_P$.

We obtain a homogeneous $A$-module of finite type. By virtue of Theorem 3.1, it is extremely
projective with respect to the category ${\cal H}$.

Now take $\e>0$. As we know, for every $n$ the operator $\s_n$ has a right inverse operator, say $\rho_n$, with the norm $<1+\e/2$. Consider the morphisms $\s:Q\to P$ and $\rho:P\to Q$, generated by the sequences $\s_n$ and $\rho_n$, respectively. Obviously, $\rho$ is a right inverse to $\s$. By virtue of Proposition 1.10(ii), $f_Q=f_P$, and we can apply to both morphisms Proposition 1.11. Therefore, since all $\s_n$ are contractive, the same is true for
$\s$, and the estimate for numbers $\|\rho_n\|$ gives $\|\rho\|\le\e/2<\e$. Thus $\s$ is a near--retraction of normed modules, and all what remains is to apply Proposition 1.1.

\bigskip
{{\bf The end of the proof of Theorem I}}

\bigskip
Now we suppose that a module $P\in{\cal H}$ is extremely projective with respect to the category ${\cal H}$

Of course, every coordinate submodule of a homogeneous module is itself homogeneous. Therefore
we can use Proposition 1.4. This immediately implies that $P$ has the property (i).

To complete the proof of Theorem I, we proceed to show that every module in ${\cal H}$, which is topologically  projective with respect to ${\cal H}$, has finite type.

In what follows, for an arbitrary $A$--module $X$ we denote, for brevity, the set $X \setminus X_{00}$ by $X^\ii$.

Suppose that for $X\in{\cal H}$ we have $X^\ii\ne\emptyset$. Our nearest aim is to construct, with the help of $X$, another module ${\bf Y}\in{\cal H}$ that will serve as the domain of a future coisometry onto $X$.

Let us choose and fix, for every $x\in X^\ii$, an arbitrary $\bar x\in\overline{X}$ and the sequence $\lm_n^x$ with the properties, indicated in Proposition 1.7. Further, for the same $x$, denote by $Y^x$ the submodule $\{\mu\bar x+b\cd\bar x; \mu\in\C, b\in A\}$ in $\overline{X}$, in other words the submodule in $\overline{X}$, {\it algebraically} generated by $\bar x$.

After this, we introduce the $A$-module
$$
{\bf Y}:=\bigoplus\{Y^x;x\in X^\ii\}\bigoplus X_{00},
$$
the {\it algebraic} direct sum of the indicated modules. In what follows, for a given ${\bf y}\in {\bf Y}$
the notation ${\bf y}=\{y^{x_k};k=1,...,m, y^*\}$ means that $y^{x_k}\in Y^{x_k}$ and $y^*\in X_{00}$ are direct summands of this element, and all other summands are zeroes. Note that ${\bf y}_n$, the $n$-th coordinate of our ${\bf y}$, can be presented as $\{y^{x_k}_n;k=1,...,m, y^*_n\}$.

\medskip
We want to equip ${\bf Y}$ with a norm. To avoid a misunderstanding, we shall denote this future norm by $\|\cd\|_{\bf Y}$, whereas the already given norm on $X$, as well as on its completion $\overline{X}$, will be denoted just by $\|\cd\|$.

We begin with the coordinate submodules of {\bf Y}. Fix, for a moment, $n\in\N$ and set, for an arbitrary ${\bf y}=\{y^{x_k};k=1,...,m, y^*\}\in {\bf Y}_n$ (where, of course, $y^{x_k}\in Y^{x_k}_n$ and $y^*\in X_n$),
$$
\|{\bf y}\|_n:=\sum_{k=1}^m\|y^{x_k}_n\|+\|y^*_n\|.
$$
Evidently, $\|\cd\|_n$ is a norm on ${\bf Y}_n$.

Make the following obvious observation: for every $n\in\N$, we have $X_n=0$ if and only if ${\bf Y}_n=0$. It follows that, in the notation of Section 1, $c_{00}^+({\bf Y}_{00})=c_{00}^+(X_{00})$. Denote by $f$ the function, associated with the module $X_{00}$.
Now, using the recipe of Proposition 1.9, we introduce, with the help of that $f$ and the norms $\|\cd\|_n$, the norm on ${\bf Y}_{00}$. Denote this norm by $\|\cd\|^0_{\bf Y}$.

Thus ${\bf Y}_{00}$ becomes a homogeneous $A$-module. Now take an arbitrary ${\bf y}\in{\bf Y}$ and consider numbers $\|{\bf P}^N\cd{\bf y}\|^0_{\bf Y}$ for all $N\in\N$. By homogeneity, they form an increasing sequence. Further, note that for every ${\bf y}\in{\bf Y}_{00}$ with the only non-zero direct summand, say $y$, we have $\|{\bf y}\|_{\bf Y}^0=\|y\|$. It follows that for every  ${\bf y}=\{y^{x_k};k=1,...,m, y^*\}\in{\bf Y}$ and $N$ we have
$$
\|{\bf P}^N\cd{\bf y}\|^0_{\bf Y}\le\sum_{k=1}^m\|{\bf P}^N\cd y^{x_k}\|+\|{\bf P}^N\cd y^*\|\le\sum_{k=1}^m\|y^{x_k}\|+\|y^*\|.
$$
Thus the sequence $\|{\bf P}^N\cd{\bf y}\|^0_{\bf Y};N=1,2,...$ converges; denote its limit by $\|{\bf y}\|_{\bf Y}$. Using the respective properties of the norm $\|\cd\|_{\bf Y}^0$, it is easy to check that $\|\cd\|_{\bf Y}$ is also a norm, this time on all ${\bf Y}$, and it makes the latter an essential homogeneous $A$-module.

Of course, if a certain ${\bf y}$ belongs to some direct summand of ${\bf Y}$, that is either to $Y^x$ for some $x\in X^\ii$, or to $X_{00}$, then $\|{\bf y}\|_{\bf Y}$ is the norm of our element in the respective submodule of $\ov{X}$. On the other hand, if ${\bf y}$ belongs, for some $n$, to ${\bf Y}_n$, that is ${\bf y}$ can be written as $\{y^{x_k};k=1,...,m, y^*\}$ with all $y^{x_k}\in(\ov{X})_n$ and $y^*\in X_n$, then we obviously have
$$
\|{\bf y}\|_{\bf Y}=\|{\bf y}\|_n=\sum_{k=1}^m\|y^{x_k}\|+\|y^*\|.\eqno(4)
$$

\medskip
Our next aim is to introduce a morphism of $A$-modules $\tau:{\bf Y}\to X$.

\smallskip
It is sufficient to define $\tau$ on direct summands of ${\bf Y}$. First, take a summand of the form $Y^x; x\in X^\ii$ and its element, say $y$. The latter, as we remember, has the form $\mu\bar{x}+b\cd\bar{x}$ for some
$\mu\in\C$ and $b\in A$. We set $\tau(y):=\mu{x}+b\cd{x}$. Since we have, of course, $\tau(y)_n=(\lm^x_n)^{-1}y_n$ (cf. Proposition 1.7), and $X$ is homogeneous, our $\tau(y)$ is uniquely defined by $y$.
This gives rise to the well defined map from $Y^x$ into $X$, which is obviously a morphism of $A$-modules. In the case of the remaining direct summand, $X_{00}$, we define $\tau$ just as the natural embedding into $X$.

Thus $\tau$ is defined. Note that for ${\bf y}\in{\bf Y}$, say ${\bf y}=\{y^{x_k};k=1,...,m, y^*\}$, we
obviously have
$$
\tau({\bf y}_n)=\sum_{k=1}^m(\lm^{x_k}_n)^{-1}y^{x_k}_n+y^*_n.\eqno(5)
$$

\medskip
{\bf Proposition 3.2}. The morphism $\tau$ is coisometric.

\smallskip
PROOF. Clearly, $\tau$ has a well defined birestriction $\tau_{00}:{\bf Y}_{00}\to X_{00}$, and $\tau_{00}$ is
generated by its coordinate submorphisms $(\tau_{00})_n:({\bf Y}_{00})_n\to (X_{00})_n$; here, of course, $({\bf Y}_{00})_n={\bf Y}_n$ and $(X_{00})_n=X_n$.

Fix, for a moment, $n\in\N$ and take ${\bf y}=\{y^{x_k}_n;k=1,...,m, y_n^*\}\in{\bf Y}_n$.
We obviously have $\tau(y^{x_k}_n)=(\lm^{x_k}_n)^{-1}y^{x_k}_n$, where, as we remember from Proposition 1.7,
$\lm^{x_k}_n\ge1$. Combining this with the equality (4), we see that the norm of $(\tau_{00})_n({\bf y})$, that is of $\tau({\bf y})$, does not exceed $\|{\bf y}\|$.

From this, by Proposition 1.11, we have $\|\tau_{00}({\bf y})\|\le\|{\bf y}\|$ for all ${\bf y}\in{\bf Y}_{00}$
and hence $\|\tau({\bf P}^N\cd{\bf y})\|\le\|{\bf P}^N\cd{\bf y}\|$ for all ${\bf y}\in{\bf Y}$ and $N\in\N$. But, by virtue of Proposition 1.5, $\tau({\bf y})=\lim_{N\to\ii}{\bf P}^N\cd\tau({\bf y})$. Besides, $\tau$ is a module morphism. From this we see that $\tau$ is contractive.

It remains to display, for a given $x\in X$ and $\e>0$, a
certain ${\bf y}\in {\bf Y}$ with $\tau({\bf y})=x$ and $\|{\bf y}\|_{\bf Y}<\|x\|+\e$. If $x\in X_{00}$, then the
copy of $x$ in the direct summand $X_{00}$ of ${\bf Y}$ already fits. If $x\in X^\ii$, we recall that $\ov{X}\in{\cal H}$ (Proposition 1.6), and therefore, by Proposition 1.5, there exists $N\in\N$ such that $\|\bar x-{\bf P}^N\cd\bar x\|<\e$.

Now consider ${\bf y}\in{\bf Y}$ with at most two non-zero direct summands, namely \\ $\bar x-{\bf P}^N\cd\bar x\in Y^x$ and ${\bf P}^N\cd x\in X_{00}$. We see that $\tau({\bf y})=
x$, and

$$
 \|{\bf y}\|_{\bf Y}\le  \|\bar x-{\bf P}^N\cd\bar x\|+\|{\bf P}^N\cd x\|<\e+\|x\|.\:\:\: \Box
$$

\medskip
{\bf Proposition 3.3.} {\it Let $X;X^\ii\ne\emptyset$, ${\bf Y}$ and $\tau$ be as before. Then
the identity morphism $\id_X$ has no bounded lifting with respect to $\tau$. In other words, there is no bounded morphism of $A$-modules $\psi: X\to {\bf Y}$, making the diagram
$$
\xymatrix@R-10pt@C+15pt{
& {\bf Y} \ar[d]^{\tau}\\
X \ar[ur]^{\psi} \ar[r]^{\id} & X }
$$
commutative.}

\smallskip
PROOF. Let $\psi$ be an (algebraic) morphism of $A$-modules, being a lifting of $\id_X$ with respect to $\tau$. Fix an arbitrary $x\in X^\ii$. Denote, for brevity, $\psi(x)$ by ${\bf y}$ and write it as ${\bf y}=\{y^{x_k};k=1,...,m, y^*\}$ (see above).
Since $\psi$ is a morphism, $\psi(x_n)$ coincides with ${\bf y}_n$ and thus can be written as
$\{y^{x_k}_n;k=1,...,m, y^*_n\}$, where $y^{x_k}_n\in (Y^{x_k})_n$ and $y^*_n\in X_n$.
Since $y^*$ is a finite sequence, the equality (4) transforms to $\|{\bf y}_n\|=\sum_{k=1}^m\|y^{x_k}_n\|$ for sufficiently large $n$. But by (5) for these $n$ we have also
$$
x_n=\tau\psi(x_n)=\tau({\bf y}_n)=\sum_{k=1}^m(\lm^{x_k}_n)^{-1}y^{x_k}_n.
$$
Therefore for all sufficiently big $n$ we have
$$
\|x_n\|\le\zeta_n\sum_{k=1}^m\|y^{x_k}_n\|=\zeta_n\|{\bf y}_n\|_{\bf Y},
$$
where we set $\zeta_n:=\max\{(\lm^{x_k}_n)^{-1}; k=1,...,m\}$. Since $\zeta_n$ tends to 0, the set of numbers
$\|\psi(x_n)\|_{\bf Y}/\|x_n\|$, taken over all $n$ with $x_n\ne0$, is not bounded.
This shows, of course, that the morphism $\psi$ is not bounded. $\Box$

\bigskip
Recall that the module {\bf Y} in the previous proposition belongs, together with $X$, to the category ${\cal H}$. Therefore this proposition shows that the last assertion of Theorem I is valid. Since every extremely projective module is certainly topologically projective, we obtain
the `only if' part of the theorem into the bargain. This completes the proof of Theorem I.

\bigskip
{\bf The proof of Theorem II}

\bigskip
{\bf `If' part.} Suppose we are given a module $X\in\ov{{\cal H}}$ with $X_n=l_1(\Lambda_n); n=1,2,\dots$. Consider its submodule $P$ of finite type with $P_n:=l_1^0(\Lambda_n)$. Of course, $P$ belongs to ${\cal H}$, and hence, by Theorem 3.1, it is extremely projective with respect to that category. But the completion of $P$ is obviously our initial $X$. Therefore, by Propositions 1.6 and 1.2 combined, $X$ is  extremely projective with respect to $\ov{{\cal H}}$.

\medskip
{\bf `Only if' part.} Now we suppose that $Q\in\ov{{\cal H}}$ is extremely projective with respect to $\ov{{\cal H}}$. Clearly, $\ov{{\cal H}}$ satisfies the condition on $\f$, formulated in Proposition 1.5, and $\f_n$, for $\f:=\ov{{\cal H}}$, is {\bf Ban}. Consequently, the mentioned proposition implies that for every $n$, the Banach space $Q_n$ is  extremely projective in {\bf Ban}. It remains to apply the Grothendieck Theorem, formulated in the Introduction.


\ed